\documentclass[a4paper,12pt]{article}
\usepackage[a4paper,margin=2.5truecm]{geometry}
\usepackage{amsmath,amsfonts,amssymb,amsthm,txfonts,graphicx,hyperref,url}
\newtheorem{definition}{Definition}
\newtheorem{lemma}{Lemma}
\newtheorem{hypothesis}{Hypothesis}
\newtheorem{openproblem}{Open Problem}
\newtheorem{theorem}{Theorem}
\newcommand{\N}{\mathbb N}
\title{A common approach to three open problems in number theory}
\author{Apoloniusz Tyszka\\
{\normalsize University of Agriculture}\\
{\normalsize Faculty of Production and Power Engineering}\\
{\normalsize Balicka 116B, 30-149 Krak\'ow, Poland}\\
{\normalsize E-mail address:} \url{rttyszka@cyf-kr.edu.pl}}
\begin{document}
\begin{sloppypar}
\date{}
\maketitle
\begin{abstract}
The following system of equations \mbox{$\{x_1 \cdot x_1=x_2$}, \mbox{$x_2 \cdot x_2=x_3$},
\mbox{$2^{\textstyle 2^{\textstyle x_1}}=x_3$}, \mbox{$x_4 \cdot x_5=x_2$}, \mbox{$x_6 \cdot x_7=x_2\}$}
has exactly one solution in \mbox{$(\N\ \setminus \{0,1\})^7$}, namely \mbox{$(2,4,16,2,2,2,2)$}.
Hypothesis~\ref{hyp1} states that if a system of equations
\mbox{${\cal S} \subseteq \{x_i \cdot x_j=x_k: i,j,k \in \{1,\ldots,7\}\} \cup$}
\mbox{$\{2^{\textstyle 2^{\textstyle x_j}}=x_k: j,k \in \{1,\ldots,7\}\}$}
has at most five equations and at most finitely many solutions in \mbox{$(\N \setminus \{0,1\})^7$},
then each such solution \mbox{$(x_1,\ldots,x_7)$} satisfies \mbox{$x_1,\ldots,x_7 \leqslant 16$}.
Hypothesis~\ref{hyp1} implies that there are infinitely many
composite numbers of the form \mbox{$2^{\textstyle 2^{\textstyle n}}+1$}.
Hypotheses~\ref{hyp2} and \ref{hyp3} are of similar kind.
Hypothesis~\ref{hyp2} implies that if the equation \mbox{$x!+1=y^2$}
has at most finitely many solutions in positive integers $x$ and $y$,
then each such solution \mbox{$(x,y)$} belongs to the set \mbox{$\{(4,5),(5,11),(7,71)\}$}.
Hypothesis~\ref{the3} implies that if the equation \mbox{$x(x+1)=y!$} has at most finitely many
solutions in positive integers $x$ and $y$, then each such solution \mbox{$(x,y)$} belongs
to the set \mbox{$\{(1,2),(2,3)\}$}. We describe \mbox{semi-algorithms} ${\tt sem}_j$ \mbox{$(j=1,2,3)$}
that never terminate. For every \mbox{$j \in \{1,2,3\}$},
if Hypothesis~$j$ is true, then ${\tt sem}_j$ endlessly prints consecutive positive
integers starting from $1$. For every \mbox{$j \in \{1,2,3\}$},
if Hypothesis~$j$ is false, then ${\tt sem}_j$ prints
a finite number (including zero) of consecutive positive integers starting from $1$.
\end{abstract}
\section{\normalsize Composite numbers of the form \mbox{$2^{\textstyle 2^{\textstyle n}}+1$}}
Let ${\cal A}$ denote the following system of equations:
\[
\Bigl\{x_i \cdot x_j=x_k:~i,j,k \in \{1,\ldots,7\}\Bigr\} \cup \Bigl\{2^{\textstyle 2^{\textstyle x_j}}=x_k:~j,k \in \{1,\ldots,7\}\Bigr\}
\]
\newpage
\noindent
The following subsystem of ${\cal A}$
\begin{center}
\includegraphics[scale=0.047]{d1.pdf}
\end{center}
\noindent
has exactly one solution in \mbox{$(\N \setminus \{0,1\})^7$}, namely \mbox{$(2,4,16,2,2,2,2)$}.
\begin{hypothesis}\label{hyp1}
If a system of equations \mbox{${\cal S} \subseteq {\cal A}$} has at most five equations and
at most finitely many solutions in \mbox{$(\N \setminus \{0,1\})^7$}, then each such solution
\mbox{$(x_1,\ldots,x_7)$} satisfies \mbox{$x_1,\ldots,x_7 \leqslant 16$}.
\end{hypothesis}
\begin{lemma}\label{lem1} (\cite[p.~109]{Tyszka1}).
For every \mbox{non-negative} integers $x$ and $y$, \mbox{$x+1=y$} if and only if
\mbox{$2^{\textstyle 2^{\textstyle x}} \cdot 2^{\textstyle 2^{\textstyle x}}=2^{\textstyle 2^{\textstyle y}}$}.
\end{lemma}
\begin{theorem}\label{the1}
Hypothesis~\ref{hyp1} implies that \mbox{$2^{\textstyle 2^{\textstyle x_1}}+1$} is composite
for infinitely many integers $x_1$ greater than~$1$.
\end{theorem}
\begin{proof}
Assume, on the contrary, that Hypothesis~\ref{hyp1} holds and \mbox{$2^{\textstyle 2^{\textstyle x_1}}+1$}
is composite for at most finitely many integers $x_1$ greater than~$1$. Then, the equation
\[
x_2 \cdot x_3=2^{\textstyle 2^{\textstyle x_1}}+1
\]
has at most finitely many solutions in \mbox{$(\N \setminus \{0,1\})^3$}.
By Lemma~\ref{lem1}, in positive integers greater than~$1$, the following subsystem of ${\cal A}$
\begin{center}
\includegraphics[scale=0.047]{d2.pdf}
\end{center}
has at most finitely many solutions in \mbox{$(\N \setminus \{0,1\})^7$} and expresses that
\[
\left\{\begin{array}{rcl}
x_2 \cdot x_3 &=& 2^{\textstyle 2^{\textstyle x_1}}+1 \\
x_4 &=& 2^{\textstyle 2^{\textstyle x_1}}+1 \\
x_5 &=& 2^{\textstyle 2^{\textstyle x_1}} \\
x_6 &=& 2^{\textstyle 2^{\textstyle 2^{\textstyle 2^{\textstyle x_1}}}} \\
x_7 &=& 2^{\textstyle 2^{\textstyle 2^{\textstyle 2^{\textstyle x_1}}+1}}
\end{array}\right.
\]
Since \mbox{$641 \cdot 6700417=2^{\textstyle 2^{\textstyle 5}}+1>16$}, we get a contradiction.
\end{proof}
\newpage
Most mathematicians believe that \mbox{$2^{\textstyle 2^{\textstyle n}}+1$} is composite for
every integer \mbox{$n \geqslant 5$}, see \mbox{\cite[p.~23]{analytic}}.
\begin{openproblem} (\mbox{\cite[p.~159]{17lectures}}).
Are there infinitely many composite numbers of the form \mbox{$2^{\textstyle 2^{\textstyle n}}+1$}?
\end{openproblem}
Primes of the form \mbox{$2^{\textstyle 2^{\textstyle n}}+1$} are called Fermat primes,
as Fermat conjectured that every integer of the form \mbox{$2^{\textstyle 2^{\textstyle n}}+1$}
is prime, see \mbox{\cite[p.~1]{17lectures}}. Fermat remarked that \mbox{$2^{\textstyle 2^{\textstyle 0}}+1=3$},
\mbox{$2^{\textstyle 2^{\textstyle 1}}+1=5$}, \mbox{$2^{\textstyle 2^{\textstyle 2}}+1=17$},
\mbox{$2^{\textstyle 2^{\textstyle 3}}+1=257$}, and \mbox{$2^{\textstyle 2^{\textstyle 4}}+1=65537$}
are all prime, see \mbox{\cite[p.~1]{17lectures}}.
\begin{openproblem} (\mbox{\cite[p.~158]{17lectures}}).
Are there infinitely many prime numbers of the form \mbox{$2^{\textstyle 2^{\textstyle n}}+1$}?
\end{openproblem}
\section{\normalsize The Brocard-Ramanujan equation \mbox{$x!+1=y^2$}}
Let ${\cal B}$ denote the following system of equations:
\[
\{x_i \cdot x_j=x_k: i,j,k \in \{1,\ldots,6\}\} \cup \{x_j!=x_k: (j,k \in \{1,\ldots,6\}) \wedge (j \neq k)\}
\]
The following subsystem of ${\cal B}$
\begin{center}
\includegraphics[scale=0.047]{d3.pdf}
\end{center}
has exactly two solutions in positive integers, namely \mbox{$(1,\ldots,1)$} and \mbox{$(2,2,4,24,24!,(24!)!)$}.
\begin{hypothesis}\label{hyp2}
If a system of equations \mbox{${\cal S} \subseteq {\cal B}$} has at most finitely many
solutions in positive integers \mbox{$x_1,\ldots,x_6$}, then each such solution
\mbox{$(x_1,\ldots,x_6)$} satisfies \mbox{$x_1,\ldots,x_6 \leqslant (24!)!$}.
\end{hypothesis}
\begin{lemma}\label{lem2}
For every positive integers $x$ and $y$, \mbox{$x! \cdot y=y!$} if and only if
\[
(x+1=y) \vee (x=y=1)
\]
\end{lemma}
\begin{theorem}\label{the2}
Hypothesis~\ref{hyp2} implies that if the equation \mbox{$x_1!+1=x_2^2$}
has at most finitely many solutions in positive integers $x_1$ and $x_2$,
then each such solution \mbox{$(x_1,x_2)$} belongs to the set \mbox{$\{(4,5),(5,11),(7,71)\}$}.
\end{theorem}
\begin{proof}
The following system of equations ${\cal B}_1$
\begin{center}
\includegraphics[scale=0.047]{d4.pdf}
\end{center}
\newpage
\noindent
is a subsystem of ${\cal B}$. By Lemma~\ref{lem2}, in positive integers, the system ${\cal B}_1$ expresses
that \mbox{$x_1=\ldots=x_6=1$} or
\[
\left\{\begin{array}{rcl}
x_1!+1&=& x_2^2 \\
x_3 &=& x_1! \\
x_4 &=& (x_1!)! \\
x_5 &=& x_1!+1 \\
x_6 &=& (x_1!+1)!
\end{array}\right.
\]
\par
If the equation \mbox{$x_1!+1=x_2^2$}
has at most finitely many solutions in positive integers $x_1$ and $x_2$, then ${\cal B}_1$
has at most finitely many solutions in positive integers \mbox{$x_1,\ldots,x_6$} and
Hypothesis~\ref{hyp2} implies that every tuple \mbox{$(x_1,\ldots,x_6)$} of positive integers
that solves ${\cal B}_1$ satisfies \mbox{$(x_1!+1)!=x_6 \leqslant (24!)!$}. Hence, \mbox{$x_1 \in \{1,\ldots,23\}$}.
If \mbox{$x_1 \in \{1,\ldots,23\}$}, then \mbox{$x_1!+1$} is a square only for \mbox{$x_1 \in \{4,5,7\}$}.
\end{proof}
\par
It is conjectured that \mbox{$x!+1$} is a square only for \mbox{$x \in \{4,5,7\}$},
see \mbox{\cite[p.~297]{Weisstein}}.
A weak form of Szpiro's conjecture implies that the equation \mbox{$x!+1=y^2$}
has only finitely many solutions in positive integers, see \cite{Overholt}.
\section{\normalsize Erd\"os' equation \mbox{$x(x+1)=y!$}}
Let ${\cal C}$ denote the following system of equations:
\[
\{x_i \cdot x_j=x_k: (i,j,k \in \{1,\ldots,6\}) \wedge (i \neq j)\} \cup
\{x_j!=x_k: (j,k \in \{1,\ldots,6\}) \wedge (j \neq k)\}
\]
The following subsystem of ${\cal C}$
\begin{center}
\includegraphics[scale=0.047]{d5.pdf}
\end{center}
has exactly three solutions in positive integers, namely \mbox{$(1,\ldots,1)$}, \mbox{$(1,1,2,2,2,2)$}, and
\mbox{$(2,2,3,6,720,720!)$}.
\begin{hypothesis}\label{hyp3}
If a system of equations \mbox{${\cal S} \subseteq {\cal C}$} has at most finitely many solutions in
positive integers \mbox{$x_1,\ldots,x_6$}, then each such solution \mbox{$(x_1,\ldots,x_6$)} satisfies
\mbox{$x_1,\ldots,x_6 \leqslant 720!$}.
\end{hypothesis}
\begin{theorem}\label{the3}
Hypothesis~\ref{the3} implies that if the equation \mbox{$x_1(x_1+1)=x_2!$} has at most finitely many
solutions in positive integers $x_1$ and $x_2$, then each such solution \mbox{$(x_1,x_2)$} belongs
to the set \mbox{$\{(1,2),(2,3)\}$}.
\end{theorem}
\newpage
\begin{proof}
The following system of equations ${\cal C}_1$
\begin{center}
\includegraphics[scale=0.047]{d6.pdf}
\end{center}
is a subsystem of ${\cal C}$. By Lemma~\ref{lem2}, in positive integers,
the system ${\cal C}_1$ expresses that \mbox{$x_1=\ldots=x_6=1$} or
\[
\left\{\begin{array}{rcl}
x_1 \cdot (x_1+1)&=& x_2! \\
x_3 &=& x_1 \cdot (x_1+1) \\
x_4 &=& x_1! \\
x_5 &=& x_1+1 \\
x_6 &=& (x_1+1)!
\end{array}\right.
\]
\par
If the equation \mbox{$x_1(x_1+1)=x_2!$} has at most finitely many solutions in positive
integers $x_1$ and~$x_2$, then ${\cal C}_1$
has at most finitely many solutions in positive integers \mbox{$x_1,\ldots,x_6$} and
Hypothesis~\ref{hyp3} implies that every tuple \mbox{$(x_1,\ldots,x_6)$} of positive integers
that solves ${\cal C}_1$ satisfies \mbox{$x_2!=x_3 \leqslant 720!$}.
Hence, \mbox{$x_2 \in \{1,\ldots,720\}$}.
If \mbox{$x_2 \in \{1,\ldots,720\}$}, then $x_2!$ is a product of two consecutive positive integers
only for \mbox{$x_2 \in \{2,3\}$} because the following {\sl MuPAD} program
\begin{quote}
\begin{verbatim}
for x2 from 1 to 720 do
x1:=round(sqrt(x2!+(1/4))-(1/2)):
if x1*(x1+1)=x2! then print(x2) end_if:
end_for:
\end{verbatim}
\end{quote}
returns $2$ and $3$.
\end{proof}
\par
The question of solving the equation \mbox{$x(x+1)=y!$} was posed by \mbox{P. Erd\"os}, see \cite{Berend}.
F. Luca proved that the $abc$ conjecture implies that the equation \mbox{$x(x+1)=y!$}
has only finitely many solutions in positive integers, see \cite{Luca}.
\section{\normalsize Hypotheses~\ref{hyp2} and \ref{hyp3} cannot be generalized to an arbitrary number of variables}
Let \mbox{$f(1)=2$}, \mbox{$f(2)=4$}, and let $f(n+1)=f(n)!$ for every integer \mbox{$n \geqslant 2$}.
Let \mbox{${\mathcal W}_1$} denote the system of equations \mbox{$\{x_1!=x_1$}.
For an integer \mbox{$n \geqslant 2$}, let \mbox{${\mathcal W}_n$} denote the following system of equations:
\begin{center}
\includegraphics[width=\hsize]{d7.pdf}
\end{center}
\newpage
For every positive integer $n$, the system \mbox{${\mathcal W}_n$}
has exactly two solutions in positive integers \mbox{$x_1,\ldots,x_n$},
namely \mbox{$(1,\ldots,1)$} and \mbox{$(f(1),\ldots,f(n))$}.
For a positive integer $n$, let \mbox{$\Psi_n$} denote the following statement: {\em if a system
of equations
\[
{\mathcal S} \subseteq \{x_i \cdot x_j=x_k: i,j,k \in \{1,\ldots,n\}\} \cup \{x_j!=x_k: j,k \in \{1,\ldots,n\}\}
\]
has at most finitely many solutions in positive integers \mbox{$x_1,\ldots,x_n$},
then each such solution \mbox{$(x_1,\ldots,x_n)$} satisfies \mbox{$x_1,\ldots,x_n \leqslant f(n)$}.}
\begin{theorem}\label{the4}
Every factorial Diophantine equation can be algorithmically transformed into
an equivalent system of equations of the forms \mbox{$x_i \cdot x_j=x_k$} and \mbox{$x_j!=x_k$}.
It means that this system of equations satisfies a modified version of Lemma 4 in \cite{Tyszka1}.
\end{theorem}
\begin{proof}
It follows from Lemmas~\mbox{2--4} in \cite{Tyszka1} and Lemma~\ref{lem2}.
\end{proof}
\par
The statement \mbox{$\forall n \in \N \setminus \{0\} ~\Psi_n$} is dubious.
By Theorem~\ref{the4}, this statement implies that there is an algorithm
which takes as input a factorial Diophantine equation and returns an integer
which is greater than the solutions in positive integers, if these solutions
form a finite set. This conclusion is strange because
properties of factorial Diophantine equations are similar to properties
of exponential Diophantine equations and a computable upper
bound on non-negative integer solutions does not exist for exponential Diophantine
equations with a finite number of solutions, see \cite{Matiyasevich}.
\section{\normalsize Equivalent forms of Hypotheses~\ref{hyp1}--\ref{hyp3}}
If \mbox{$k \in [10^{19},10^{20}-1] \cap \N$}, then there are uniquely determined
\mbox{non-negative} integers \mbox{$a(0),\ldots,a(19) \in \{0,\ldots,9\}$} such that
\[
\Bigl(a(19) \geqslant 1\Bigr) \wedge \Bigr(k=a(19) \cdot 10^{19}+a(18) \cdot 10^{18}+\ldots+a(1) \cdot 10^1+a(0) \cdot 10^0\Bigr)
\]
\begin{definition}\label{def1}
For an integer \mbox{$k \in [10^{19},10^{20}-1]$}, ${\mathcal S}_k$ stands for the smallest system of
equations~${\mathcal S}$ satisfying conditions {\tt (1)} and {\tt (2)}.
\vskip 0.2truecm
\noindent
{\tt (1)} If \mbox{$i \in \{0,4,8,16\}$} and \mbox{$a(i) \in \{0,1,2,3,4\}$},
then the equation \mbox{$x_{a(i+1)} \cdot x_{a(i+2)}=x_{a(i+3)}$} belongs to ${\mathcal S}$
when it belongs to ${\mathcal A}$.
\vskip 0.2truecm
\noindent
{\tt (2)} If \mbox{$i \in \{0,4,8,16\}$} and \mbox{$a(i) \in \{5,6,7,8,9\}$},
then the equation \mbox{$2^{\textstyle 2^{\textstyle x_{a(i+1)}}}=x_{a(i+2)}$} belongs to~${\mathcal S}$
when it belongs to ${\mathcal A}$.
\end{definition}
\begin{lemma}\label{newlemma1}
$\{{\mathcal S}_k: k \in [10^{19},10^{20}-1] \cap \N\}=
\{{\mathcal S}: ({\mathcal S} \subseteq {\mathcal A}) \wedge ({\rm card}({\mathcal S}) \leqslant 5)\}$.
\end{lemma}
\begin{proof}
It follows from the equality \mbox{$5 \cdot 4=20$}.
\end{proof}
\newpage
For a positive integer $n$, let $p_n$ denote the \mbox{$n$-th} prime number.
\begin{theorem}\label{newtheorem1}
The following \mbox{semi-algorithm} ${\tt sem}_1$ never terminates.
\begin{center}
\includegraphics{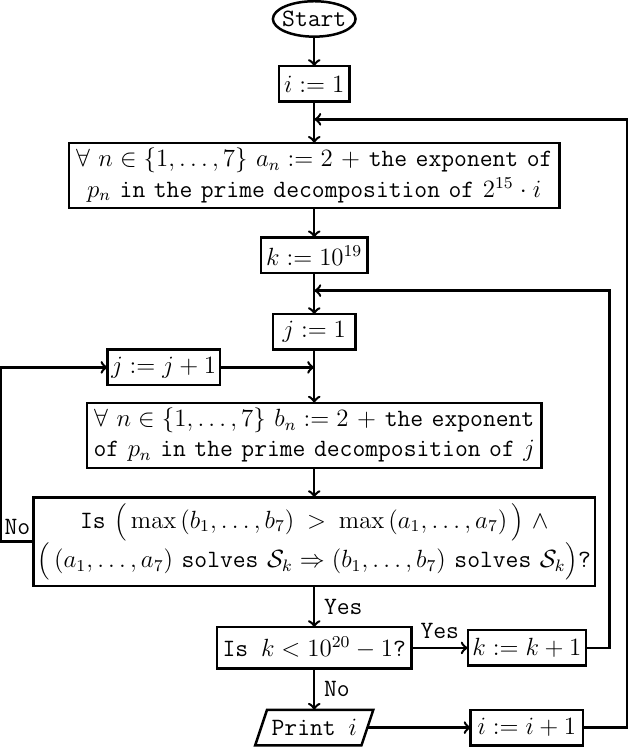}
\end{center}
\noindent
If Hypothesis~\ref{hyp1} is true, then ${\tt sem}_1$ endlessly prints consecutive positive
integers starting from $1$. If Hypothesis~\ref{hyp1} is false, then ${\tt sem}_1$ prints
a finite number (including zero) of consecutive positive integers starting from $1$.
\end{theorem}
\begin{proof}
It follows from Lemma~\ref{newlemma1}.
\end{proof}
\newpage
\begin{theorem}\label{newtheorem2}
The following \mbox{semi-algorithm} ${\tt sem}_2$ never terminates.
\begin{center}
\includegraphics{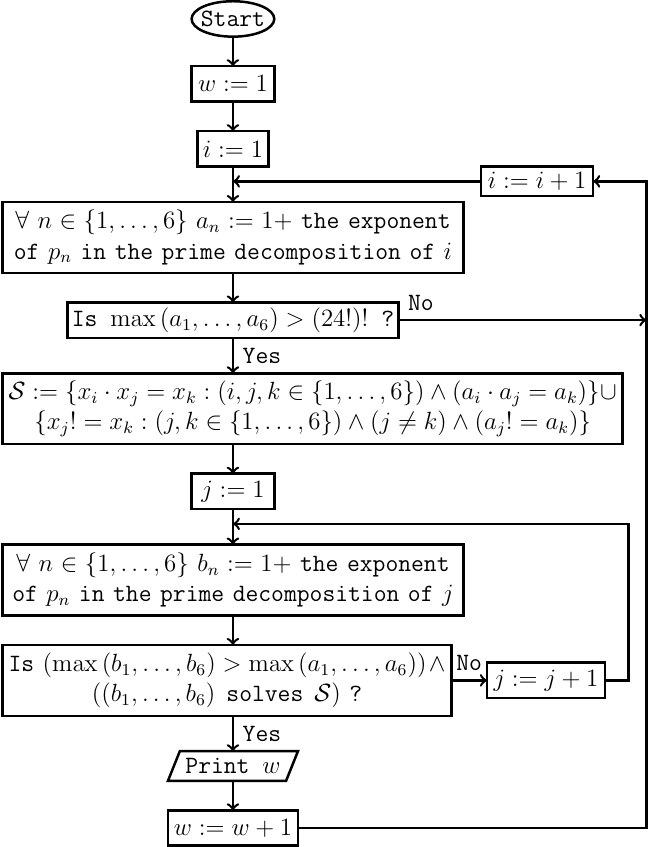}
\end{center}
\noindent
If Hypothesis~\ref{hyp2} is true, then ${\tt sem}_2$ endlessly prints consecutive positive
integers starting from $1$. If Hypothesis~\ref{hyp2} is false, then ${\tt sem}_2$ prints
a finite number (including zero) of consecutive positive integers starting from $1$.
\end{theorem}
\newpage
\begin{theorem}\label{newtheorem3}
The following \mbox{semi-algorithm} ${\tt sem}_3$ never terminates.
\begin{center}
\includegraphics{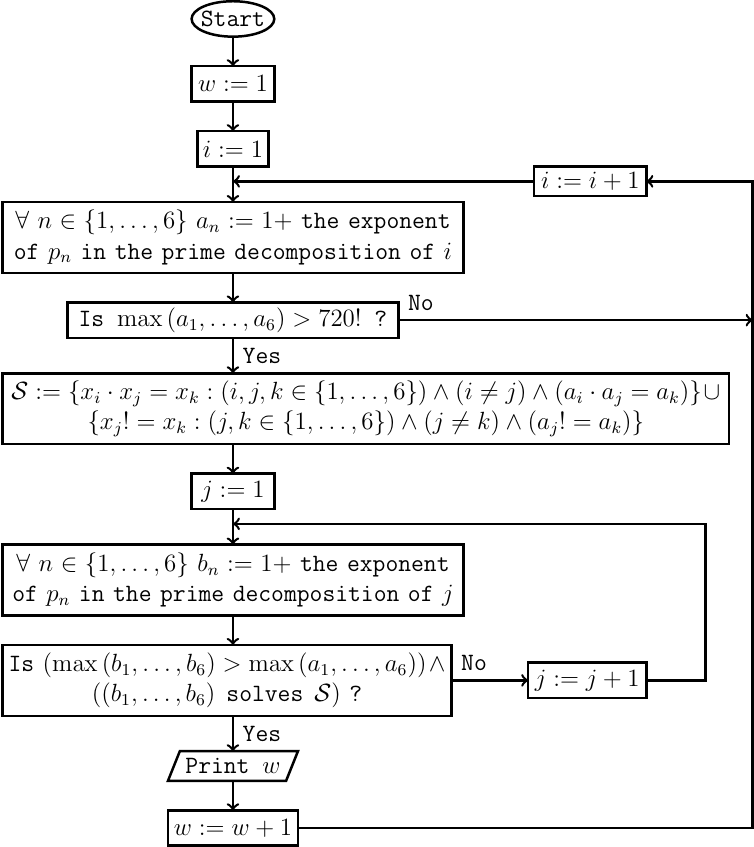}
\end{center}
If Hypothesis~\ref{hyp3} is true, then ${\tt sem}_3$ endlessly prints consecutive positive
integers starting from $1$. If Hypothesis~\ref{hyp3} is false, then ${\tt sem}_3$ prints
a finite number (including zero) of consecutive positive integers starting from $1$.
\end{theorem}

\end{sloppypar}
\end{document}